\documentclass[11pt, leqno]{amsart}

\usepackage[mathscr]{eucal}
\usepackage{amsmath}
\usepackage{amsfonts}
\usepackage{amsthm}
\usepackage{color}
\usepackage{xypic}

\theoremstyle{plain}
\newtheorem{theorem}{Theorem}[section]
\newtheorem{proposition}{Proposition}[section]
\newtheorem{remark}{Remark}

\usepackage{a4wide}
\usepackage{amssymb}
\usepackage{graphicx}

\numberwithin{equation}{section}

\title[Periodic magnetic curves]{Periodic magnetic curves in elliptic Sasakian space forms}

\author[J,~Inoguchi]{Jun-ichi Inoguchi}
\address[J.~Inoguchi] {Department of Mathematical Sciences,
Yamagata University, \newline
Yamagata, 990-8560,
Japan}
\email{inoguchi@sci.kj.yamagata-u.ac.jp}

\author[M.~I.~Munteanu]{Marian Ioan Munteanu}
\address[M.~I.~Munteanu]
{University 'Al. I. Cuza' of Iasi, 
Faculty of Mathematics, Bd. Carol I, no.~11,
700506 Iasi, Romania}
\email{marian.ioan.munteanu@gmail.com}

\date{\today}

\begin{document}

\begin{abstract}
It is an interesting question whether a given equation of motion has a periodic solution or not, and
in the positive case to describe them.
We investigate periodic magnetic curves in elliptic Sasakian space forms and
we obtain a quantization principle for periodic magnetic flowlines on Berger spheres.
We give a criterion for periodicity of magnetic curves on the unit sphere ${\mathbb{S}}^3$.
\end{abstract}



\maketitle

\section{Introduction}
It has been a long-standing problem 
to find closed trajectories in dynamical systems 
on manifolds. The existence of closed trajectories 
turns out to be subtle and closely related to 
topology and geometric structures of the manifolds.
    


For example, the following fundamental existence theorem of closed 
geodesics is well known. 
 
For any compact Riemannian manifold $(M,g)$, 
in each element in the fundamental group 
$\pi_{1}(M)$, there exists a closed geodesic which attains 
the minimum of the energy functional in the homotopy class.  

One of the interesting generalizations of the notion of geodesic is 
that of magnetic curve.

Magnetic curves represent, in physics, 
the trajectories of the charged particles moving on a 
Riemannian manifold under the magnetic forces. 
Geodesics are magnetic curves with free magnetic field.
Moreover, magnetic curves are characterized as 
critical points of the Landau-Hall functional. 

In any $3$-dimensional Riemannian manifold 
$(M,g)$, magnetic fields of nonzero constant length 
are in one to one correspondence to almost contact structure  
compatible to the metric $g$ (\textit{cf.} [3]). 
This fact motivates us to study magnetic curves in almost contact metric $3$-manifolds (with 
closed 
fundamental $2$-form). 

As we have mentioned above, 
magnetic curves are derived from the variational problem 
of the Landau-Hall functional. From variational point of view,
the existence of global potential for magnetic fields is a 
natural assumption. 
Based on these observations, we study magnetic curves in 
$3$-dimensional contact metric manifolds, especially in Sasakian manifolds. 

In \cite{Tau07}, Taubes proved the generalized Weinstein conjecture in dimension 3, namely, 
{\em on compact, orientable, contact $3$-manifold the Reeb vector field $\xi$ has at least one closed
integral curve. 
}
In conjunction with this problem it is important to explore the existence of periodic magnetic trajectories
of the contact magnetic field defined by $\xi$ in Sasakian manifolds. We know (see \cite{DIMN}) that the study of
magnetic curves in Sasakian space forms of arbitrary dimension reduces to their investigation in dimension $3$.

 In 2009 Cabrerizo et al. \cite{CFG09}
have been looked for periodic orbits of the contact magnetic field on the unit sphere ${\mathbb{S}}^3$.
See also \cite{Bar08}.
In this paper we study closed magnetic curves in arbitrary elliptic Sasakian space forms.

\section{Preliminaries}

\subsection{Magnetic curves} The motion of the charged particles in a Riemannian manifold under the action of 
the magnetic fields are known as magnetic curves. More precisely, a \textit{magnetic field} $F$ on a 
Riemannian manifold $(M, g)$ is a closed $2$-form $F$ and the \textit{Lorentz force} associated to $F$ 
is a tensor field $\phi$ of type $(1,1)$ such that
\begin{equation}\label{eqLforce}
F(X,Y) = g(\phi X,Y), \ \  X,Y \in \mathfrak{X}(M).
\end{equation}
A curve $\gamma$ on $M$ that satisfies the \textit{Lorentz equation} 
\begin{equation}
\label{eqL}
\nabla_{\dot\gamma}\dot\gamma=\phi(\dot\gamma),
\end{equation}
is called \textit{magnetic trajectory} of $F$ or simply a \emph{magnetic curve}.
Here $\nabla$ denotes the Levi Civita connection associated to the metric $g$.
Lorentz equation generalizes the equation of geodesics under arc length parametrization, 
namely, $\nabla_{\dot\gamma}\dot\gamma=0$. 
A magnetic field $F$ is said to be {\em uniform} if $\nabla F = 0$.  

It is well-known that the magnetic trajectories have constant speed. When the magnetic curve $\gamma(s)$ is arc length parametrized, 
it is called a \textit{normal magnetic curve}. 

The dimension 3 is rather special, since it allows us to identify 2-forms with vector fields via the Hodge $\star$ operator 
and the volume form $dv_g$ of the (oriented) manifold. In this way, magnetic fields may be identified with divergence free 
vector fields by 
\begin{equation*}
\label{mag_FV}
F_{V}=\iota_{V}dv_{g}.
\end{equation*}
Magnetic fields $F$ corresponding to Killing vector fields are usually known as
\textit{Killing magnetic fields}. 
Their trajectories, called \textit{Killing magnetic curves}, are of great importance
since they are related to the Kirchhoff elastic rods. See e.g. \cite{BCFR07, BarrosRomero}.

\subsection{Sasakian manifolds} 
A {\em $(\varphi, \xi, \eta)$ structure} on a manifold $M$ is defined by a field $\varphi$ of endomorphisms of 
tangent spaces, a vector field $\xi$ and a 1-form $\eta$ satisfying
\begin{equation*}
\eta(\xi)=1,\ \  \varphi^{2}=-\mathrm{I}+\eta\otimes\xi, \ \ \varphi\xi=0,\ \ \eta \circ \varphi = 0.
\end{equation*}
If $(M, \varphi, \xi, \eta)$ admits a compatible Riemannian metric $g$, namely
\begin{equation*}
g(\varphi{X},\varphi{Y})=g(X,Y)-\eta(X)\eta(Y),
{\rm \ for\  all\  } X, Y\in \mathfrak{X}(M),
\end{equation*}
 then $M$ is said to have an 
{\em almost contact metric structure}, and $(M, \varphi, \xi, \eta, g)$ is called an {\em almost contact metric manifold}. 
Consequently, we have that $\xi$ is unitary
and $\eta(X) = g(\xi, X)$, for any $X\in\mathfrak{X}(M)$.

We define a 2-form $\Omega$ on $(M, \varphi, \xi, \eta, g)$ by 
\begin{equation}\label{2form}
\Omega(X,Y) = g(\varphi X, Y),
{\rm \ for\  all\  } X, Y\in \mathfrak{X}(M),
\end{equation}
called {\em the fundamental $2$-form} of the almost contact metric structure $(\varphi, \xi, \eta, g)$.

If $\Omega = d\eta$, 
then  $(M, \varphi, \xi, \eta, g)$ is called a {\em contact metric manifold}. Here $d\eta$ is defined by 
$
d\eta(X,Y) = \frac{1}{2}\big(X\eta(Y) - Y\eta(X)-\eta([X,Y])\big),
$
for any $X, Y\in \mathfrak{X}(M)$.
On a contact metric manifold $M$, the 1-form $\eta$ is a contact form (see Blair's book \cite{Blair}). 
The vector field $\xi$ is called the \textit{Reeb vector field} of $M$ and it is characterized by $\iota_\xi \eta =1$ and $\iota_\xi d\eta = 0$. 
Here $\iota$ denotes the interior product. 
In analytical mechanics $\xi$ is traditionally called the \textit{characteristic vector field} of $M$.

An almost contact metric manifold $M$ is said to be {\em normal} if the normality tensor  
\linebreak
$
\displaystyle
S= N_\varphi(X, Y)+2d\eta(X,Y)\xi
$ vanishes, where $N_\varphi$ is the \textit{Nijenhuis torsion} of $\varphi$ defined by
$
\displaystyle
N_\varphi(X, Y) = 
[\varphi X,\varphi Y ] + \varphi^2[X, Y ] 
-\varphi[\varphi X, Y ] -\varphi[X, \varphi Y ],
$
for any $X, Y \in \mathfrak{X}(M)$.

A {\em Sasakian manifold} is defined as a normal contact metric manifold. 
Denoting by $\nabla$ the Levi Civita connection associated to $g$, 
the Sasakian manifold $(M, \varphi, \xi, \eta, g)$ is characterized by 
\begin{equation*}
\label{Sasaki1}
(\nabla_X \varphi)Y = - g(X, Y)\xi + \eta(Y)X,
\ {\rm \ for\ any\ }\ X,Y\in\mathfrak{X}(M).
\end{equation*}
As a consequence, we have
\begin{equation}\label{Sasaki2}
\nabla_X \xi = \varphi X, \ \forall X\in \mathfrak{X}(M).
\end{equation}

A contact metric structure $(\varphi, \xi, \eta, g)$ is called {\em $K$-contact} if $\xi$ is a Killing vector field. 
Due to \eqref{Sasaki2} and the fact that $\varphi$ is skew-symmetric, it follows that a Sasakian manifold is $K$-contact. 
The converse is not true in general. Yet, a 3-dimensional manifold is Sasakian if and only if it is $K$-contact.

A plane section $\Pi$ at $p\in M^{2n+1}$ is called a $\varphi$-section if it is invariant under $\varphi_p$. The sectional curvature 
$k(\Pi)$ of a $\varphi$-section is called the {\em $\varphi$-sectional curvature} of $M^{2n+1}$ at $p$. A Sasakian manifold 
$(M^{2n+1}, \varphi, \xi, \eta, g)$ is said to be a {\em Sasakian space form} if it has constant $\varphi$-sectional curvature.

Take a positive constant $a$ and define a new Sasakian structure $(\varphi,\hat{\xi},\hat{\eta},\hat{g})$ on $M$ by
$$
\hat{\xi}:=\frac{1}{a}\xi,\ \
\hat{\eta}:=a\eta, \  \
\hat{g}:=ag+a(a-1)\eta\otimes\eta.
$$
This structure is called a $D$-\textit{homothetic deformation} of 
$(\varphi,\xi,\eta,g)$.
In particular, if $M(c)$ is a Sasakian space form, then deforming the structure we obtain also a Sasakian space form $M(\hat{c})$, where 
$\hat{c} = \frac{c+3}{a}-3$. For every value of $c$ there exists Sasakian space forms, as follows: 
the elliptic Sasakian space forms, know also as the Berger spheres if $c>-3$, 
the Heisenberg space $\mathbb{R}^{2n+1}(-3)$, if $c=-3$, and $B^{2n}\times\mathbb{R}$ when $c<-3$. 
See also \cite[Theorem 7.15]{Blair}.
Notice that the case $c>-3$ includes the standard unit sphere $\mathbb{S}^{2n+1}(1)$.

\subsection{Magnetic curves in Sasakian manifolds}

Let $(M, \varphi, \xi, \eta, g)$ be a contact metric manifold and let $\Omega$ be the fundamental 2-form defined by \eqref{2form}. 
Since $\Omega=d\eta$ on a contact metric manifold, $\Omega$ is a closed 2-form, thus 
we can define a magnetic field on $M$ by
\begin{equation}\label{contactmagfieldF}
F_q(X,Y) = q\Omega(X,Y), 
\end{equation}
where $X,Y\in\mathfrak{X}(M)$ and $q$ is a real constant.
We call $F_q$ the {\em contact magnetic field} with the {\em strength} $q$.
Notice that if $q=0$, then the contact magnetic field vanishes identically and the magnetic curves are the geodesics of $M$. 
In the sequel we assume $q\neq 0$.

The {\em Lorentz force} $\phi_q$ associated to the contact magnetic field $F_q$ may be easily determined combining 
\eqref{2form} and \eqref{eqLforce}, namely
\begin{equation}\label{Lforceq}
\phi_q = q\varphi, 
\end{equation} 
where $\varphi$ is the field of endomorphisms of the contact metric structure. 

In this setting, the Lorentz equation \eqref{eqL} can be written as
\begin{equation}
\label{eqLcontact}
\nabla_{\dot\gamma} \dot\gamma = q\varphi\dot\gamma,
\end{equation}
where $\gamma:I\subseteq\mathbb{R}\to M^{2n+1}$ is a smooth curve parametrized by arc length.
The solutions of \eqref{eqLcontact} are called {\em normal magnetic curves} or {\em trajectories} for $F_q$.

Let $F=F_V$ be a magnetic field with 
corresponding divergence free vector field $V$ on a 
$3$-dimensional Riemannian manifold $(M,g)$.

Then one can check that the Lorentz force $\phi$ of $F$ satisfies 
([4], page 7):
$$
\phi^{2}X=-g(V,V)X+g(V,X)V.
$$ 
Moreover, we have
$$
V^{\flat}(\phi{X})=g(\phi{X},V)=
g(V\times X,V)=dv_{g}(V,X,V)=0,
$$
$$
g(\phi{X},\phi{Y})
=g(V,V)g(X,Y)-V^{\flat}(X)V^{\flat}(Y).
$$
Here $V^\flat$ denotes the $1$-form metrically dual to $V$.
Thus if $V$ has \textit{constant} length $q>0$, 
then 
$(\varphi,\xi,\eta)=(\phi/q,V/q,V^{\flat}/q)$ defines a 
$(\varphi,\xi,\eta)$-structure compatible to the Riemannian structure $g$ 
satisfying $\mathrm{div}\xi=0$. The magnetic field $F$ is represented 
by $F=q\Omega$ as a contact magnetic field on the resulting 
almost contact metric manifold $(M,\varphi,\xi,\eta,g)$.

\section{Elliptic Sasakian space forms}
\label{section1}

In this section, we shall give an explicit matrix group model of a
simply connected elliptic Sasakian space form $\mathscr{M}^3(c)$.

As is well known, the unit $3$-sphere
$(\mathbb{S}^3;\eta_1,\xi_1,\varphi_1,g_1)$ is identified with the special
unitary group $G=\mathrm{SU}(2)$ with bi-invariant metric. In this
section we shall give a $\mathrm{SU}(2)$-model of
$\mathscr{M}^3(c)$.

Let us denote the Lie algebra of $G$ by $\mathfrak{g}$.
The bi-invariant metric $g_1$ of constant curvature $1$
on $G$ is induced by the following inner
product $\langle \cdot , \cdot \rangle_1$
on $\mathfrak{g}$:
$$
\langle X,Y \rangle_1=-\frac{1}{2}\ \mathrm{tr} (XY),\ \  X,Y \in \mathfrak{g}.
$$

Take a quaternionic basis of $\mathfrak{g}$:
$$
{\mathbf i}=
\left (
\begin{array}{cc}
0 & \sqrt{-1} \\
\sqrt{-1} & 0
\end{array}
\right ),
\ \
{\mathbf j}=
\left (
\begin{array}{cc}
0 & -1 \\
1 & 0
\end{array}
\right ),
\ \
{\mathbf k}=
\left (
\begin{array}{cc}
\sqrt{-1} & 0 \\
0 & -\sqrt{-1}
\end{array}
\right ).
$$
By using this basis, the group $\mathrm{SU}(2)$ is described as
$$
\mathrm{SU}(2)=
\left \{
\
\left (
\begin{array}{cc}
x_0+\sqrt{-1}~x_3 &
-x_2+\sqrt{-1}~x_1
\\
x_2+\sqrt{-1}~x_1
&
x_0-\sqrt{-1}~x_3
\end{array}
\right )
\ \ \biggr \vert
\ \
x_0^2+x_1^2
+x_2^2+x_3^2=1
\right \}.
$$
In the spinor representation of the Euclidean $3$-space ${\mathbb{E}}^3$,
we identify ${\mathbb{R}}^3$ with ${\mathfrak{g}}={\mathfrak{su}}(2)$
via the correspondence
$$
(x_1,x_2,x_3)\longleftrightarrow
x_1\mathbf{i}+x_2\mathbf{j}+x_3\mathbf{k}=
\left(
\begin{array}{rr}
\sqrt{-1} ~x_3&  -x_2+\sqrt{-1}~x_1\\[2mm]
x_2+\sqrt{-1}~ x_1& -\sqrt{-1}~x_3
\end{array}
\right).
$$
Denote the left translated vector fields of $\{{\mathbf i},
{\mathbf j}, {\mathbf k}\}$ by $\{E_1,E_2,E_3\}$. Then a left
invariant Sasakian structure of $G$ is given by
$$
\xi_1:=E_3, \ \ \eta_1=g_1(E_3,\cdot ),
$$
$$
\varphi_1(E_1)=-E_2,\ \ \varphi_1(E_2)=E_1,\ \
\varphi_1(E_3)=0.
$$
Note that the commutation relations of
$\{E_1,E_2,E_3\}$ are
$$
[E_1,E_2]=2E_3,\ \
[E_2,E_3]=2E_1,\ \
[E_3,E_1]=2E_2.
$$
The Lie group $G$ acts isometrically on the Lie algebra
$\mathfrak{g}$ by the $\mathrm{Ad}$-action.
$$
\mathrm{Ad}:G \times \mathfrak{g}
\rightarrow \mathfrak{g};
\ \mathrm{Ad}(a)X=aXa^{-1},
\ \ a \in G,\ X \in \mathfrak{g}.
$$
The $\mathrm{Ad}$-orbit of $\mathbf{k}/2$ is a sphere of radius
$1/2$ in the Euclidean $3$-space $\mathbb{E}^3=\mathfrak{g}$. The
$\mathrm{Ad}$-action of $G$ on $S^2(1/2)$ is isometric and
transitive. The isotropy subgroup of $G$ at $\mathbf{k}/2$ is
$$
\mathrm{U}(1)=\left \{
a_t=
\left (
\begin{array}{cc}
e^{\sqrt{-1}t}
& 0
\\ 0 &
e^{-\sqrt{-1}t}
\end{array}
\right )
\ \biggr \vert
\ t \in {\mathbb R}
\right \}.
$$
Hence $\mathbb{S}^2(1/2)$ is represented by $\mathrm{SU}(2)/\mathrm{U}(1)$
as a Riemannian symmetric space. The natural projection
$$
\pi_1:{\mathbb{S}}^3 \rightarrow \mathbb{S}^2(1/2),
\ \ \pi_1(a)=\mathrm{Ad}(a)(\mathbf{k}/2)
$$
is a Riemannian submersion and
defines a principal $\mathrm{U}(1)$-bundle
over $\mathbb{S}^2(1/2)$.

Since the Sasakian structure $(\eta_1,\xi_1,\varphi_1,g_1)$ is left
invariant, its $D$-homothetic deformation is also left invariant.
Hence the elliptic contact space form $\mathscr{M}^3(c)$ is
identified with $\mathrm{SU}(2)$ endowed with the left invariant contact
metric structure:
$$
\eta:=\alpha \eta_1,\ \
\xi:=\xi_1/\alpha,\ \
\varphi:=\varphi_1,
$$
$$
g(X,Y)
=\alpha
g_{1}(X,Y)+
\alpha(\alpha-1)
\eta_1(X)\eta_1(Y),\ \
c=4/\alpha-3,
$$
where $\alpha$ is a positive real number.
The Reeb vector field $\xi$ generates a one parameter group of
transformations on $\mathscr{M}^3(c)$. Since $\xi$ is a Killing
vector field, this transformation group acts isometrically on $G$.
The transformation group generated by $\xi$ is identified with the
following Lie subgroup $K$ of $G$:
$$
K=\left \{
\left (
\begin{array}{cc}
e^{\frac{\sqrt{-1}}{\alpha}t}
& 0
\\ 0 &
e^{-\frac{\sqrt{-1}}{\alpha}t}
\end{array}
\right )
\ \Biggr \vert
\ t \in {\mathbb R}
\right \}
\cong
\mathrm{U}(1).
$$
Furthermore, the action of the transformation group generated by
$\xi$ corresponds to the natural right action of $K$ on $G$:
$$
G \times K \rightarrow G;\
(a,k) \mapsto ak.
$$
By using the curvature formula due to O'Neill, one can
see that the orbit space $G/K$ is a $2$-sphere of radius
$\sqrt{\alpha}/2$, namely the  constant curvature $(c+3)$-sphere. 
The Riemannian metric $g$ is not
only $G$-left invariant but also $K$-right invariant. Hence
$G\times K$ acts isometrically on $G$. The elliptic contact space
form $\mathscr{M}^3(c)$ is represented by $G\times K/K=G$ as a
naturally reductive homogeneous space. For $c\not=1$,
$\mathscr{M}^3(c)$ has $4$-dimensional isometry group.

In particular, $g$ is $G$-bi-invariant if and only if $c=1$. In this
case $\mathscr{M}^3(1)$ is represented by $G\times G/G$ as a
Riemannian symmetric space. Note that $\mathscr{M}^3(1)$ has
$6$-dimensional isometry group.

\vspace{0.2cm}

Now we shall take an orthonormal frame field $\{e_1,e_2,e_3\}$ of
$\mathscr{M}^3(c)$ by
$$
e_1:=\frac{1}{\sqrt{\alpha}}E_1,\ \
e_2:=\frac{1}{\sqrt{\alpha}}E_2,\ \
e_3:=\frac{1}{\alpha}\xi_1.
$$
Then the commutation relations of this basis are
$$
[e_1,e_2]=2e_3,\ \
[e_2,e_3]=\frac{c+3}{2}e_1,\ \
[e_3,e_1]=\frac{c+3}{2}e_2.
$$

The Levi-Civita connection $\nabla$ of $(\mathscr{M}^3(c),g)$ is
described by
$$
\nabla_{e_1}e_1=0,\ \
\nabla_{e_1}e_2=e_3,\ \
\nabla_{e_1}e_3=-e_2,
$$
$$
\nabla_{e_2}e_1=-e_3,\ \
\nabla_{e_2}e_2=0,\ \
\nabla_{e_2}e_3=e_1,
$$
$$
\nabla_{e_3}e_1=\frac{c+1}{2}e_2,\ \
\nabla_{e_3}e_2=-\frac{c+1}{2}e_1,\ \
\nabla_{e_3}e_3=0.
$$

The Riemannian curvature
tensor field $R$ of $(\mathscr{M}^{3}(c),g,\nabla)$
is described by
\begin{equation*}
R_{1212}=c,\ \ R_{1313}=R_{2323}=1,
\end{equation*}
and the sectional curvatures are:
\begin{equation*}
K_{12}=c,
\  K_{13}=K_{23}=1.
\end{equation*}
 The Ricci tensor $Ric$ and the
scalar curvature $scal$ are computed to be
\begin{equation*}
Ric_{11}=Ric_{22}
=c+1, \ Ric_{33}=2,
\ scal=2(c+2).
\end{equation*}

\section{Magnetic trajectories in 3-dimensional elliptic Sasakian space forms}
Let us consider a normal magnetic curve with respect to 
the magnetic field $F_q=q\Phi$, where $q\in{\mathbb{R}}$ is the strength and $\Phi$ is the fundamental $2$-form 
defined by $\Phi(X,Y)=g(\varphi X,Y)$. Thus, the Lorentz force is $\phi=q\varphi$.

Let us expand the tangent vector field $T(s)=\gamma^{\prime}(s)$ as
$T(s)=T_{1}(s)e_{1}+T_{2}(s)e_{2}+T_{3}(s)e_3$. 
Then the tension field $\tau(\gamma)$ is 
$$
\tau(\gamma)=\nabla_{\gamma^\prime}
\gamma^{\prime}
=\left\{T_{1}^{\prime}-\frac{1}{2}(c-1)T_2T_3\right\}e_1
+\left\{T_{2}^{\prime}+\frac{1}{2}(c-1)T_1T_3\right\}e_2+
T_{3}^{\prime}e_3.
$$ 
From the third equation we have 
$\cos\theta:=T_{3}$ is constant.

Hence, the magnetic equation 
$\tau(\gamma)=q\varphi{T}$ is
\begin{align*}
T_{1}^{\prime}(s)&=~\tilde{q}\>T_{2}(s), \\
T_{2}^{\prime}(s)&=-\tilde{q}\>T_{1}(s),
\end{align*}
where we put
$\tilde{q}= q+\frac{1}{2}(c-1)\cos \theta.
$

Integrating,  we obtain 
\begin{align*}
T_{1}(s)&= \cos (\tilde{q}s)\>T_{1}(0)+\sin (\tilde{q}s)\>T_{2}(0),\\
T_{2}(s)&= \sin (\tilde{q}s)\>T_{1}(0)-\cos (\tilde{q}s)\>T_{2}(0),\\
T_{3}(s)&=\cos\theta.
\end{align*}
One can see that $\kappa=|q|\sin\theta$.

If we consider $\alpha=1$, then $c=1$ and $\tilde{q}=q$.

\begin{remark} \rm
If $\gamma(s)$ is a normal magnetic curve in $\mathscr{M}^{3}(1)$, for a certain strength $q$,
then, after an affine change of the parameter, $\gamma$ becomes a normal magnetic curve in $\mathscr{M}^{3}(c)$ 
with a different strength $\hat q$. See for details \cite{DIMN}.
\end{remark}

\bigskip

As we have seen in Section~\ref{section1} we have a circle bundle
$$
\xymatrix{
{\mathrm{U}}(1) \ar[r] & \mathscr{M}^{3}(c) \ar[d] ^\pi \\ & {\mathbb{S}}^2(r)}
$$
where the projection $\pi$ is defined by $\pi(a)={\mathrm{Ad}}(a)(r\mathbf{k})$.
If ${\mathbb{S}}^2(r)\subset {\mathfrak{su}}(2)\equiv{\mathbb{E}}^3$ is endowed with the induced metric from $\langle~,~\rangle_1$
then $\pi$ becomes a Riemannian submersion from $( \mathscr{M}^{3}(c),g)$ onto ${\mathbb{S}}^2(r)$ if and only if
$r=\frac{\sqrt{\alpha}}{2}$.

The correspondence $a\mapsto \ker\eta_a$, $a\in  \mathscr{M}^{3}(c)$ defines a connection in this principal ${\mathrm{U}}(1)$ bundle, hence
if $A\in\ker\eta_a$ then $|\pi_{*,a}(A)|=|A|$. Take $\eta_a$ as the connection form (the standard choice).

Consider a regular curve $\beta:{\mathbb{R}}\longrightarrow {\mathbb{S}}^2(r)\subset \big({\mathfrak{su}}(2),\langle~,~\rangle_1\big)$,
$u\mapsto \beta(u)$. As usual, $\beta$ will be parametrized by the arc length and let $\hat\beta$ be a horizontal
lift of $\beta$. This means that $\pi(\hat\beta(u))=\beta(u)$ for all $u\in{\mathbb{R}}$ and
$\langle\hat\beta(u)^*\hat\beta'(u),{\mathbf{k}}\rangle=0$.
Here $\hat\beta(u)^*=\overline{\hat\beta(u)^T}$, and $a^T$ denotes the transpose of the matrix $a$.
The inner product on ${\mathfrak{su}}(2)$, the tangent space of $\mathscr{M}^{3}(c)$ at the identity, is 
$\langle X,Y\rangle=\alpha\langle X,Y\rangle_1+\alpha(\alpha-1)\langle X,{\mathbf{k}}\rangle_1\langle Y,{\mathbf{k}}\rangle_1$.

The complete lift of $\beta$, namely $\pi^{-1}(\beta)$ is a flat surface in ${\mathrm{SU}}(2)$ and it is usually called 
\emph{the Hopf tube} over $\beta$. 
In \cite{Kit88} we find a method for constructing all flat tori in ${\mathbb{S}}^3$. See also \cite{Kit03}.
Denote $\pi^{-1}(\beta)$ by $H_\beta$. It can be naturally parametrized by
$$
F:{\mathbb{R}}\times{\mathbb{R}}\longrightarrow {\mathrm{SU}}(2),\quad
(t,u)\longmapsto F(t,u)=\hat\beta(u) a_t.
$$

In fact we can prove that the metric $g$ restricted to $H_\beta$ can be expressed as
\begin{equation*}
\label{metric_H_beta}
g_{H_\beta}=\alpha^2dt^2+du^2.
\end{equation*}
Hence we have proved the following:
\begin{proposition}
If $\beta$ is a curve on ${\mathbb{S}}^2(r)$ of length $L$, then the corresponding Hopf tube $H_\beta$
is isometric to ${\mathbb{S}}^1(\alpha)\times[0,L]$, where ${\mathbb{S}}^1(\alpha)$ is the unit circle endowed with the metric
$\alpha^2dt^2$. Moreover, its mean curvature in $ \mathscr{M}^{3}(c)$ is $H=\frac12\kappa_\beta$, where
$\kappa_\beta$ is the geodesic curvature of $\beta$ in ${\mathbb{S}}^2(r)$.
\end{proposition}

If $\beta$ is a closed curve, i.e. $\beta(u+L)=\beta(u)$ for all $u\in{\mathbb{R}}$, then the relation
$F(t,u)=\hat\beta(u)a_t$ defines a covering of the $(t,u)$ plane onto the immersed torus in $ \mathscr{M}^{3}(c)$,
called the \emph{Hopf torus} corresponding to $\beta$. One can easily see that, if $\beta$ is a great circle in ${\mathbb{S}}^2(r)$,
then the Hopf torus $H_\beta$ is minimal in $ \mathscr{M}^{3}(c)$.

\begin{proposition}
The magnetic curve $\gamma$ is a geodesic of the Hopf tube $H_\beta$.
\end{proposition}
\proof
We know that 
$$
\nabla_{\gamma'}\gamma'=q(T_2e_1-T_1e_2),
$$
which is orthogonal both to $\xi$ and  $\gamma'$.
It follows that $\nabla_{\gamma'}\gamma'$ is collinear to the normal vector to $H_\beta$.
Hence, $\gamma$ is a geodesic on $H_\beta$.

\endproof

\section{Periodic magnetic trajectories on $\mathscr{M}^{3}(c)$}

In this section we show that the set of all periodic magnetic trajectories on $\mathscr{M}^{3}(c)$ can be quantized in the
set of rational numbers. 

First of all, we find a relation between the intrinsic geometry of the magnetic curve and those of its projection on the sphere 
${\mathbb{S}}^2(r)$. We state the following result.

\begin{proposition}
If $\gamma$ is a magnetic curve on  $\mathscr{M}^{3}(c)$, then the projection curve $\beta$ is a circle on ${\mathbb{S}}^2(r)$
of geodesic curvature
$$
\kappa_\beta=\frac{\kappa^2+\tau^2-1}{\kappa}~,
$$
where $\kappa$ and $\tau$ are (constant) curvature and torsion of $\gamma$ in $ \mathscr{M}^{3}(c)$.
\end{proposition}
\proof
For any $a\in  \mathscr{M}^{3}(c)\equiv {\mathrm{SU}}(2)$ we have
$$
{e_1}_{\big|_a}=\frac 1{\sqrt{\alpha}}~a{\mathbf{i}},
\quad
{e_2}_{\big|_a}=\frac 1{\sqrt{\alpha}}~a{\mathbf{j}},
\quad
{e_3}_{\big|_a}=\frac 1{\alpha}~a{\mathbf{k}}.
$$
Let $\gamma$ be a magnetic curve on $\mathscr{M}^{3}(c)$ such that 
$\dot \gamma=T_1e_1+T_2e_2+\cos\theta e_3$, as in the previous Section.
If $\beta=\pi(\gamma)$, then
$$
\beta'=\pi_{*,\gamma}(\dot\gamma)=r(\dot\gamma{\mathbf{k}}\gamma^*+\gamma{\mathbf{k}}\dot\gamma^*)=
     \frac{2r}{\sqrt{\alpha}}\big(-T_1\gamma{\mathbf{j}}\gamma^*+T_2\gamma{\mathbf{i}}\gamma^*\big).
$$

Since $r=\sqrt{\alpha}/2$ we obtain
\begin{equation*}
\label{betap}
\beta'=\gamma(-T_1{\mathbf{j}}+T_2{\mathbf{i}})\gamma^*.
\end{equation*}

In the same manner we can compute
$$
\beta''=\gamma V(s) \gamma^*,
$$
where
$$
V(s)=\big(T_2'+\frac2\alpha\cos\theta T_1\big){\mathbf{i}}+\big(-T_1'+\frac2\alpha\cos\theta T_2\big){\mathbf{j}}
	-\frac2{\sqrt{\alpha}}\sin^2\theta{\mathbf{k}}.
$$

If $\stackrel{2}{\nabla}$ denotes the Levi Civita connection on ${\mathbb{S}}^2(r)$ we have
$$
\beta''=\stackrel{2}{\nabla}_{\beta'}\beta'-\frac{\sin^2\theta}{r^2}~\beta.
$$
Since $T_1'=\tilde q~T_2$ and $T_2'=-\tilde q~ T_1$, we get
\begin{equation}
\label{eq:6}
\stackrel{2}{\nabla}_{\beta'}\beta'=(-q+2\cos\theta)\gamma\big(T_1{\mathbf{i}}+T_2{\mathbf{j}}\big)\gamma^*.
\end{equation}

If $\nu_1$ denotes the first (unit) normal of $\gamma$, we have
$$
\nu_1=\frac{\varepsilon}{\sqrt{\alpha}\sin\theta}~\gamma\big(T_2{\mathbf{i}}-T_1{\mathbf{j}}\big),
$$
where $\varepsilon={\rm sgn}(q)$.

Thus, the unit normal to $\beta$ can be considered
$$
\nu_\beta:=\pi_{*,\gamma}(\nu_1)=-\frac{\varepsilon}{\sin\theta}\gamma\big(T_1{\mathbf{i}}+T_2{\mathbf{j}}\big)\gamma^*.
$$
Combing this with \eqref{eq:6}, we obtain
\begin{equation*}
\label{eq:nabla2}
\stackrel{2}{\nabla}_{\beta'}\beta'=\varepsilon\sin\theta(q-2\cos\theta)\nu_\beta.
\end{equation*}
Therefore, the signed geodesic curvature of $\beta$ is
\begin{equation*}
\label{kappab}
\kappa_\beta=\frac{q-2\cos\theta}{\varepsilon\sin\theta}~.
\end{equation*}
Concerning the magnetic curve $\gamma$, we know 
\cite{DIMN}
that its curvature (the first curvature) is 
$\kappa=\varepsilon q\sin\theta$ and its torsion (the signed second curvature) is
$\tau=q\cos\theta-1$.

Since $q^2=\kappa^2+(\tau+1)^2$, the conclusion follows immediately.

\endproof

\medskip

Consider now a periodic magnetic trajectory $\gamma$ on $\mathscr{M}^{3}(c)$, i.e.
$\gamma(s+L)=\gamma(s)$ for all $s\in{\mathbb{R}}$. Then, the projection $\beta$ is a 
circle on the sphere ${\mathbb{S}}^2(r)$.

Conversely, take $\beta$ a closed curve on ${\mathbb{S}}^2(r)$ of length $L$ enclosing an
oriented area $A$, $A\in[-2\pi r^2,2\pi r^2]$. Let $\hat\beta$ be a horizontal lift of $\beta$.
Since $\hat\beta(0)$ and $\hat\beta(L)$ belongs to the same fibre, we have
$$
\hat\beta(L)=\hat\beta(0)a_\delta,\quad \delta\in{\mathbb{R}}.
$$
Usually, $\delta$ is called the holonomy number of the canonical principal connection
defined in the circle bundle. See e.g. \cite{Pink85, ABG00}. For the classical Hopf fibration, $\hat\beta$ closes up if and only if
there exists a positive integer $m$ such that, after $m$ consecutive liftings of $\beta$ we get
$\hat\beta(mL)=\hat\beta(0)$ and hence $m\delta=2\pi p$, with $p\in{\mathbb{Z}}$.

We have seen that the Hopf torus $H_\beta$ is isometric to ${\mathbb{S}}^1(\alpha)\times[0,L]$.
The isometry type depends either on the length $L$ of $\beta$, or on the area $A$ enclosed by $\beta$
on ${\mathbb{S}}^2(r)$, namely $A=\displaystyle\int\limits_D\!\!\!\int d{\mathcal{A}}$,
where $d{\mathcal{A}}$ is the area element of ${\mathbb{S}}^2(r)$, and $D$ is a domain
on ${\mathbb{S}}^2(r)$ such that $\partial D=\beta$.
Notice that $H_2({\mathbb{S}}^2(r))={\mathbb{Z}}$ and ${\mathrm{area}}({\mathbb{S}}^2(r))=4\pi r^2$.

We have:
\begin{proposition}
Let $\beta$ be a closed curve on ${\mathbb{S}}^2(r)$ of length $L$ enclosing an oriented area $A$.
Then, the corresponding Hopf torus $H_\beta$ is isometric to ${\mathbb{R}}/\Gamma$, where the lattice
$\Gamma$ is generated by the vectors $(2\pi\alpha,0)$ and $\left(\frac{A}{2r^2},L\right)$.
\end{proposition}
\proof
With previous notations we have that $\hat\beta(L)=\hat\beta(0)a_\delta$, for some $\delta$.
It is clear (due to the form of the metric on $H_\beta$) that the group of deck transformations for the
covering mentioned above is generated by the translations $(2\pi\alpha,0)$ and $(\delta,L)$.

We have to find $\delta$.
It is known that $\delta=\displaystyle\int\limits_D\!\!\!\int\Omega$, where $\Omega$ is the curvature 2-form of the connection in the
principal circle bundle $\pi:\mathscr{M}^{3}(c)\longrightarrow {\mathbb{S}}^2(r)$. It should be a multiple of the area element on 
${\mathbb{S}}^2(r)$, and one can compute
$\Omega=\frac1{2r^2}d{\mathcal{A}}$. 
Hence $\delta=\frac{A}{2r^2}$.

\endproof

If $\beta$ is the projection of a magnetic curve on $\mathscr{M}^{3}(c)$, then it is a circle on ${\mathbb{S}}^2(r)$.
Denote by $R$ its radius, $R\leq r$. We have
$$
\kappa_\beta=\frac{\sqrt{r^2-R^2}}{rR}\ ,\quad L=2\pi R,\quad A=2\pi r(r-\sqrt{r^2-R^2}).
$$
Since $\gamma$ is a periodic (closed) geodesic on the Hopf torus $H_\beta$, it corresponds to a segment in 
${\mathbb{R}}^2$ (with identified ends). This segment is in fact the diagonal of a parallelogram constructed by
taking $m$ vectors in the fibre, hence $m$ times $(2\pi\alpha,0)$ and $n$ vectors in the horizontal direction, 
i.e. $n$ times $(\frac{A}{2r^2},L)$, $n\in{\mathbb{N}}$. Thus, the direction of the magnetic trajectory $\gamma$ is given by
$$
\left(2\pi m\alpha+n\pi\Big(1-\sqrt{1-\frac{R^2}{r^2}}~\Big),2\pi nR
\right).
$$

If we put $\sigma=\cot\theta$ (here $\theta$ is the contact angle of the curve $\gamma$) and call this quantity the
\emph{slope} of $\gamma$, we have
$$
\sigma=\frac{2 m\alpha+n\Big(1-\sqrt{1-\frac{R^2}{r^2}}~\Big)}{2 nR}\ .
$$
Hence we get
$$
R\sigma+\frac12 \sqrt{1-\frac{R^2}{r^2}} =\frac mn\alpha+\frac12\ .
$$

We can state the following result.

\begin{theorem}
The set of all periodic magnetic curves of arbitrary strength on the Sasakian space form $\mathscr{M}^{3}(c)$ can be
quantized in the set of rational numbers.
\end{theorem}
\proof
Taking into account that $r=\frac{\sqrt{\alpha}}2$ we conclude that the slope of a periodic
magnetic trajectory corresponding to the circle $\beta$ of radius $R$ on the sphere ${\mathbb{S}}^2(r)$
satisfies the following quantization principle
\begin{equation}
R\sigma+\frac12\sqrt{1-\frac{4R^2}{\alpha}} ={\mathbb{Q}}\alpha+\frac12\ .
\end{equation}

\endproof

\section{Periodic magnetic curves on  ${\mathbb{S}}^3$}

We know that magnetic curves in Sasakian manifolds are helices of order 3 (see \cite{DIMN}). 
In particular, for the unit sphere ${\mathbb{S}}^3$ we have a model helix given by

\begin{equation}
\label{modelhelix_S3}
\big(\cos\psi\cos(a s), \cos\psi\sin(a s), \sin\psi\cos(b s),\sin\psi\sin(b s)\big),
\end{equation}

where $s$ is the arc-length parameter and $a,b$ and $\psi$ satisfy

$$
a^2\cos^2\psi+b^2\sin^2\psi=1.
$$

See e.g. \cite{Glu88, Tam04}. This helix lies on a torus whose mean curvature is constant $H=\cot(2\psi)$.
In fact, every helix with both curvature and torsion different from zero is congruent to one of these helices.

It is not difficult to show that the helix \eqref{modelhelix_S3} is periodic if and only if
$$
a=1/{\sqrt{p^2\sin^2\psi+\cos^2\psi}}\ ,\ b=p  a,
$$
where $p$ is a rational number. Moreover, all these helices project (via Hopf fibration) to (small) horizontal circles on
${\mathbb{S}}^2(\frac12)$.

In the following we give an example of a periodic magnetic trajectory on the unit sphere ${\mathbb{S}}^3$ and 
we draw the picture of its stereographic projection on ${\mathbb{R}}^3$.

Consider the curve $\gamma:I\longrightarrow{\mathbb{S}}^3$ defined by

\begin{equation}
\label{Ikawa_curve}
\gamma(s)=\big(x_1(s),x_2(s),x_3(s),x_4(s)\big),
\end{equation}

with
$$
\begin{array}{l}
x_1(s)=\cos(\frac{s}{2})\cos(\omega s)-\frac{1}{\omega}\big(\cos\theta-\frac12\big)\sin(\frac{s}{2})\sin(\omega s)\\[2mm]
x_2(s)=\sin(\frac{s}{2})\cos(\omega s)-\frac{1}{\omega}\big(\cos\theta-\frac12\big)\cos(\frac{s}{2})\sin(\omega s)\\[2mm]
x_3(s)=\frac{\sin\theta}{\omega}\cos(\frac{s}2)\sin(\omega s)\\[2mm]
x_4(s)=\frac{\sin\theta}{\omega}\sin(\frac{s}2)\sin(\omega s),
\end{array}
$$
where $\omega=\sqrt{\frac54-\cos\theta}$ and $\theta$ is constant. Then, $\gamma$ is a normal magnetic curve
corresponding to the contact magnetic field of ${\mathbb{S}}^3$ and  with strength 1, for which $\theta$ expresses its contact angle.
According to \cite{Ika04}, the curve $\gamma$ is periodic if and only if $\omega$ is a rational number.
In fact the periodicity condition obtained by Ikawa in \cite{Ika04} is nothing but the quantization result obtained by Cabrerizo et al.
in \cite{CFG09}.


In the following we take $\theta=\arccos\frac{29}{37}$, equivalently to $\omega=\frac23$.
Consider the stereographic projection of the sphere from its North pole. 
Then the image of $\gamma$ is drawn in the next figure.

\begin{figure}[hbtp]
\label{fig:gamma04}
\begin{center}
  \includegraphics[width=90mm]{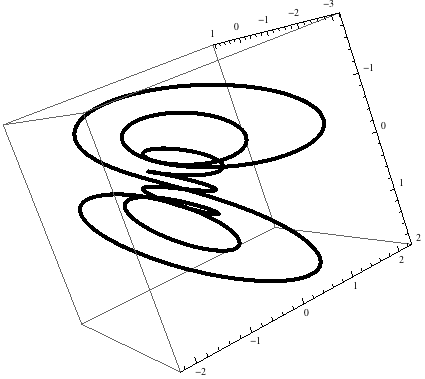}
\end{center}
\caption{$\cos\theta=29/37$}
\end{figure}

We know that $\gamma$ lies on a Hopf tube in ${\mathbb{S}}^3$. 
In the following picture we plot the image of this tube under the stereographic projection we have mentioned before.
\begin{figure}[hbtp]
\label{fig:tube04}
\begin{center}
  \includegraphics[width=60mm]{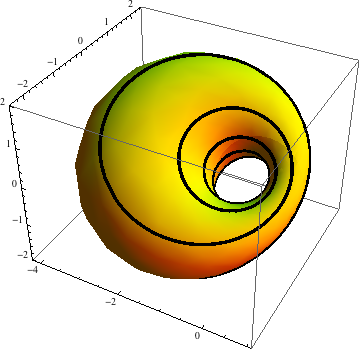}
  \quad
  \includegraphics[width=60mm]{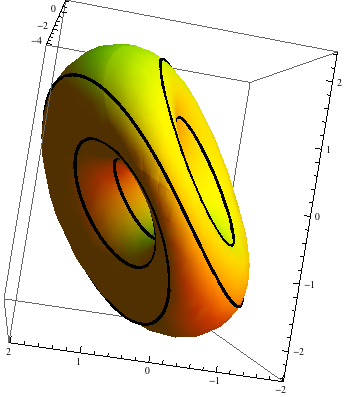}
\end{center}
\caption{The curve and the tube after stereographic projection; two different view-points}
\end{figure}

\pagebreak

We conclude this paper with a criterion of periodicity for magnetic curves in ${\mathbb{S}}^3$.
As a consequence of Theorem 3.13 in \cite{Ika04} we may state the following.
\begin{theorem}
Let $\gamma$ be a normal magnetic curve on the unit sphere ${\mathbb{S}}^3$. Then $\gamma$ is periodic if and
only if 
$$
\frac{q}{\sqrt{q^2-4q\cos\theta+4}}\in{\mathbb{Q}},
$$
where $q$ is the strength and $\theta$ is the constant contact angle of $\gamma$.
\end{theorem}

\bigskip

{\bf Acknowledgements.} 
This work was partially supported by Kakenhi 24540063 (Japan) and
by CNCS-UEFISCDI (Romania) PN-II-RU-TE-2011-3-0017/2011-2014.
The second author wishes to thank Yamagata University for warm hospitality he received
during his visit in September 2013.


\begin{thebibliography}{99}

\bibitem{ABG00} J.~Arroyo, M.~Barros, and O.~J.~Garay,
	\emph{Some examples of critical points for the total mean curvature functional},
	Proc. Edinburgh Math. Soc., {\bf 43} (2000), 587--603.

\bibitem{Bar08} M.~Barros,
	\emph{Simple geometrical models with applications in Physics},
	CP1002, Curvature and Variational Modelling in Physics and Biophysics,
	Eds. O.J. Garay, E. Garcia-Rio and R. V\'azquez-Lorenzo, AIP 2008, 71--113.
	
\bibitem{BCFR07} M.~Barros, J.~L.~Cabrerizo, M.~Fern\'andez, A.~Romero,
    {\em Magnetic vortex filament flows},
    J. Math. Phys. {\bf 48} (2007) 8, 082904:1--27.

\bibitem{BFLM} M.~Barros, A.~Ferr\'andez, P.~Lucas and M.~A.~Mero\~no,
	{\em Helicoidal filaments in the 3-sphere},
	preprint 1995, unpublished.

\bibitem{BarrosRomero} M.~Barros, A.~Romero,
    {\em Magnetic vortices},
    EPL {\bf 77} (2007), 34002:1--5.
    
\bibitem{Blair} D.~E.~Blair,  
			\emph{Riemannian Geometry of Contact and Symplectic Manifolds}, 
			Progress in Math. 203, 2002, Birkh\"auser, Boston-Basel-Berlin.
    
\bibitem{CFG09} J.~L.~Cabrerizo,  M.~Fern\'andez, and J.~S.~G\'omez,
    \emph{The contact magnetic flow in $3D$ Sasakian manifolds},
    J. Phys. A: Math. Theor. {\bf 42} (2009), 19, 195201:1--10.

\bibitem{DIMN} S.~L.~Dru\c t\u a-Romaniuc, J.~Inoguchi, M.~I.~Munteanu, and A.~I.~Nistor,
	\emph{Magnetic curves in Sasakian and cosymplectic manifolds},
	preprint.

\bibitem{Glu88} H.~Gluck,
	\emph{Geodesics in the unit tangent bundle of a round sphere},
	L'Enseignement Math\'ematique, {\bf 34} (1988) 233--246.

\bibitem{Ika04} O.~Ikawa, 
	\emph{Motion of charged particles in homogeneous K\"ahler and homogeneous Sasakian manifolds}, 
	Far East J. Math. Sci. (FJMS), 14 (3) (2004) 283--302.

\bibitem{Kit88} Y.~Kitagawa,
	\emph{Periodicity of the asimptotic curves on flat tori in ${\mathbb{S}}^3$},
	J. Math. Soc. Japan, {\bf 40} (1988) 3, 457--476.

\bibitem{Kit03} Y.~Kitagawa,
	\emph{Deformable flat tori in ${\mathbb{S}}^3$ with constant mean curvature},
	Osaka J. Math, {\bf 43} (2003), 103--119.
	
\bibitem{Pink85} U.~Pinkall,
	\emph{Hopf tori in ${\mathbb{S}}^3$},
	Invent. Math., {\bf 81} (1985), 379--386.	
	
\bibitem{Tam04} M.~Tamura,
	{\em Surfaces which contain helical geodesics in the 3-sphere},
	Mem.  Fac. Sci. Eng. Shimane Univ., {\bf 37} (2004) 59--65. 
	
\bibitem{Tau07} C.~H.~Taubes, 
	{\em The Seiberg-Witten equations and the Weinstein conjecture},
	Geom. Topol. {\bf 11} (2007), 2117--2202.

\end{thebibliography}
\end{document}